\documentclass[reqno,12pt,centertags]{article}
\usepackage{amsmath,amsthm,amscd,amssymb,latexsym,upref}
\date{\today}
\usepackage{amsmath}
\usepackage{amssymb}
\usepackage{amsfonts}
\usepackage[retainorgcmds]{IEEEtrantools} 
\usepackage{verbatim}
\usepackage{amsthm}
\usepackage{xcolor,graphicx}
\usepackage{mathrsfs}
\usepackage[bookmarks]{hyperref}

\input epsf
\usepackage{epsfig}
\usepackage[T2A,OT1]{fontenc}
\usepackage[ot2enc]{inputenc}% This one is the main point, but it does not work
\usepackage[russian,english]{babel}

\usepackage{epic,eepicemu}

\DeclareMathOperator{\sn}{sn}
\DeclareMathOperator{\cn}{cn}
\DeclareMathOperator{\dn}{dn}
\newcommand{\ud}{\,\mathrm{d}}
\newcommand{\ii}{\mathrm{i}}

\title{A Generalization of a Minimal Problem of Korkin-Zolotarev kind}
\author{N.~Akhiezer (\emph{Kharkov})}

\begin{document}

\maketitle

\begin{abstract}
This is an English translation of the paper in which N. I. Akhiezer discovered his famous orthogonal polynomials on two intervals  
in a connection with a generalization of the Korkin-Zolotarev (Korkine-Zolotaref) problem (see the small commentary by P. Yuditskii, attached to the current publication). Translated from German by F. Puchhammer.
\end{abstract}

\S 1. If one intends to approximate a continuous function $f(x)$ defined on the interval $[-1,1]$ by polynomials, one has to define 
the \emph{distance} between two continuous functions on $[-1,1]$ first, and then try to obtain the coefficients of the polynomial
$P_n(x)$, such that the distance between $f(x)$ and $P_n(x)$ becomes as small as possible.

Usually, either the quantity
\begin{equation*}
 \max_{-1\leq x\leq1}\left|f(x)-P_n(x)\right|
\end{equation*}
or the quantity
\begin{equation*}
 \int_{-1}^{1}\left[f(x)-P_n(x)\right]^2 p(x)\ud x,
\end{equation*}
where $p(x)$ denotes a fixed non-negative function on $[-1,1]$ (weight function), is chosen as a distance.

However, also the distance
\begin{equation}
 \label{eqn:one}
\int_{-1}^{1}\left|f(x)-P_n(x)\right|p(x)\ud x,
\end{equation}
which has a simple geometric meaning for $p(x)=1$, has been considered in literature several times 
(\cite{KorZolSur1873, BerSur1927, Mapkos,StiJet1876,TchSur1859}).
It is a remarkable fact that, for a large class of functions $f(x)$, which, among others, also includes $f(x)=x^{n+1}$,
the polynomial $P_n(x)$, which turns the integral in (\ref{eqn:one}) into a minimum, is just the interpolation polynomial for
$f(x)$ w.r.t. a fixed node system independent of $f(x)$.
This explains the meaning of the following problem: \emph{Among all the polynomials} $P_n(x)=x^n+p_1x^{n-1}+\cdots+p_n$ 
\emph{ find
the one for which the integral}
\begin{equation}
 \label{eqn:two}
\int_{-1}^{1}|P_n(x)|\cdot p(x)\cdot\ud x
\end{equation}
\emph{attains its smallest value.}
 This problem was first solved (for $p(x)=1$) by A. Korkine and G. Zolotaref in their paper: 
``Sur un certain minimum''.

To our knowledge, the case of an arbitrary weight function has not yet been studied.

Also, the asymptotic properties (for $n\to\infty$) of the polynomial $P_n(x)$, which turns the integral in (\ref{eqn:two}) into
a minimum, are not known to us. There is only a conjecture by S.~N.~Bernstein (\cite[p.~135]{BerSurIX}), saying that
\begin{equation*}
 \lim_{n\to\infty} \frac{H_n[t(x)]}{L_n[t(x)]}=2,
\end{equation*}
where 
\begin{equation*}
 H_n[t(x)] = \min\int_{-1}^{1}\left|x^n+p_1x^{n-1}+\cdots+p_n\right|\frac{t(x)}{\sqrt{1-x^2}}\ud x
\end{equation*}
and 
\begin{equation*}
 L_n[t(x)]=\min\left\{\max_{-1\leq x\leq 1} \left|t(x)(x^n+p_1x^{n-1}+\cdots+p_n)\right| \right\},
\end{equation*}
and where $t(x)$ is a function which is subject to certain conditions and otherwise arbitrary (trigonometric weight according
to S.~Bernstein).

Within this paper we intend to present a full solution to the above-mentioned problem in the case where
\begin{equation*}
 p(x)=\begin{cases}1 & \text{for } -1\leq x\leq\alpha \text{ and } \beta\leq x\leq1,\\ 
 0 & \text{for } \alpha<x<\beta,\end{cases}
\end{equation*}
where  $\alpha$ and $\beta$ ($-1<\alpha<\beta<1$) are two given numbers.

This case seems to be particularly interesting, since the weight function vanishes identically on a subinterval of the domain 
of integration, while, in corresponding investigations on orthogonal polynomials, this is usually only permitted on a set of mass zero
(S. Bernstein, G. Szeg\"o).

The solution is expressed in terms of elliptic functions and it is remarkable that elliptic functions suffice to obtain the 
exact solution for any $n$ and any $\alpha,~\beta$, while, in corresponding problems on obtaining the best uniform approximation
on two intervals, elliptic functions only provide the asymptotic relations (for $n\to\infty$) and only yield the exact relations for
special positions of $\alpha$ and $\beta$ (\cite{AchUeb1933}).

\S 2. So, our problem is the following:
\emph{Among all polynomials}
\begin{equation*}
 f(x)=x^n+p_1 x^{n-1}+\cdots+p_n
\end{equation*}
\emph{find the one for which the expression}
\begin{equation}
\label{eqn:three}
 I[f]=\int_{-1}^{\alpha}|f(x)|\ud x+\int_{\beta}^{1}|f(x)|\ud x=\int_{E}|f(x)|\ud x
\end{equation}
\emph{attains its smallest value, where} $\alpha,~\beta~(-1<\alpha<\beta<1)$ \emph{are given numbers and $E$ denotes the set
composed of the intervals}\footnote{The paper \cite{Pos1880} by K. Poss\' e induced me to write this one. A special problem is
posed and solved therein, which is related to the one considered here.}
\begin{equation*}
 [-1,\alpha],\quad[\beta,1].
\end{equation*}

It is clear and does not require a proof that there is at least one solution to this problem.

A repetition of the ideas of Korkine and Zolotaref shows that all zeros of the polynomial we are looking for are simple
and lie in the interval $[-1,1]$. Moreover, such a  polynomial  cannot have more than one zero between $\alpha$ and $\beta$;
indeed, assuming
\begin{equation*}
 f(x)=(x-\xi)(x-\eta)\varphi(x),
\end{equation*}
where $\alpha<\xi<\eta<\beta$, we could form the polynomial
\begin{equation*}
 f_1(x)=f(x)-h\varphi(x),
\end{equation*}
where
\begin{equation*}
 0<h<\mu=\min_{x\in E}(x-\xi)(x-\eta),
\end{equation*}
and would evidently obtain
\begin{equation*}
 I[f_1]=I[f-h\varphi]=I[f]-hI[\varphi]<I[f].
\end{equation*}

However, a zero of the polynomial under consideration lying in the interval $(\alpha,\beta)$ can also only occur in specific cases. Indeed,
let
\begin{equation*}
 f(x)=(x-\xi)\varphi(x)\qquad (\alpha<\xi<\beta),
\end{equation*}
then we have
\begin{IEEEeqnarray*}{rcl}
 I[f] &=& \int_{-1}^{\alpha}(\xi-x)|\varphi(x)|\ud x+\int_{\beta}^{1}(x-\xi)|\varphi(x)|\ud x\\
&=&
\xi\left\{ \int_{-1}^{\alpha}|\varphi(x)|\ud x - \int_{\beta}^{1}|\varphi(x)|\ud x \right\}
-\int_{-1}^{\alpha}x|\varphi(x)|\ud x + \int_{\beta}^{1}x|\varphi(x)|\ud x.
\end{IEEEeqnarray*}
Now, if the expression in braces were different from zero, then, by an appropriate variation of $\xi$, we could make
the quantity $I[f]$ smaller, i.e. $f$ could not be a solution to the problem; thus, the identity
\begin{equation*}
 \int_{-1}^{\alpha}|\varphi(x)|\ud x = \int_{\beta}^{1}|\varphi(x)|\ud x
\end{equation*}
has to hold in this case (cf. \cite{Pos1880}).

Hence, if the sought polynomial $f(x)$ has a zero $\xi$ between $\alpha$ and $\beta$, then we have the equation
\begin{equation*}
  \int_{-1}^{\alpha}\frac{|f(x)|}{\xi-x}\ud x = \int_{\beta}^{1}\frac{|f(x)|}{x-\xi}\ud x,
\end{equation*}
and, evidently, in this case, our problem has infinitely many solutions of the form
\begin{equation*}
 F(x)=\frac{x-\vartheta}{x-\xi} f(x),
\end{equation*}
where $\vartheta$ is an arbitrary number satisfying
\begin{equation*}
 \alpha\leq\vartheta\leq\beta.
\end{equation*}
The structure of these solutions shows that we can replace the zero $\xi$ by $\alpha$ or $\beta$, without altering the quantity 
$I[f]$ and we may therefore confine ourselves to such polynomials, whose zeros are simple and lie in the intervals
$[-1,\alpha]$, $[\beta,1]$.

\S 3. Let
\begin{equation*}
 \xi_1<\xi_2<\ldots<\xi_n\quad(-1\leq\xi_1,~\xi_n\leq1)
\end{equation*}
be all zeros of the polynomial being sought and let further
\begin{equation*}
 \xi_k\leq\alpha,\quad\beta\leq\xi_{k+1},
\end{equation*}
where $k$ is any of the numbers $0,1,2,\ldots,n$ and we interpret the case where $k=0$ (resp. $k=n$) as $\xi_0=-1$ (resp. $\xi_{n+1}=1$).

The integral (\ref{eqn:three}) assumes the form
\begin{IEEEeqnarray*}{rCl}
 I[f] &=& (-1)^{n}\left\{\int_{-1}^{\xi_1}f\ud x-\int_{\xi_1}^{\xi_2}f\ud x+\cdots + (-1)^{k}\int_{\xi_k}^{\alpha}f\ud x+\right.\\
  &&\left. +(-1)^k\int_{\beta}^{\xi_{k+1}}f\ud x+\cdots+(-1)^{n}\int_{\xi_n}^1f\ud x\right\}.
\end{IEEEeqnarray*}
Thus, if $p_1,~p_2,\ldots,~p_n$ denote the coefficients of the polynomial $f(x)$, then the usual conditions
\begin{equation*}
 \frac{\partial}{\partial p_i}I[f]=0\quad(i=1,2,\ldots,n)
\end{equation*}
for the extremum assume the form of the following equations, which may serve to determine the zeros $\xi_1,\xi_2,\ldots,\xi_n$:
\begin{IEEEeqnarray}{c}
 \int_{-1}^{\xi_1}x^i\ud x-\int_{\xi_1}^{\xi_2}x^i\ud x+\cdots+(-1)^{k}\int_{\xi_k}^{\alpha}x^i\ud x+ \label{eqn:four}\\
\IEEEnonumber
 +(-1)^k \int_{\beta}^{\xi_{k+1}}x^i\ud x+(-1)^{k+1}\int_{\xi_{k+1}}^{\xi_{k+2}}x^i\ud x+\cdots+(-1)^{n}\int_{\xi_n}^{1}x^i\ud x=0\\
(i=0,1,2,\ldots,n-1).\IEEEnonumber
\end{IEEEeqnarray}
Solving these equations is a very difficult task. To this end, let us notice that, due to Tch\'ebyschef and Korkine-Zolotaref \cite{Pos1880},
 by virtue of the equations in (\ref{eqn:four}),
 the relation
\begin{IEEEeqnarray*}{c}
 \int_{-1}^{\xi_1}\frac{\ud z}{x-z}-\int_{\xi_1}^{\xi_2}\frac{\ud z}{x-z}+\cdots+(-1)^{k}\int_{\xi_k}^{\alpha}\frac{\ud z}{x-z}+\\
 +(-1)^k \int_{\beta}^{\xi_{k+1}}\frac{\ud z}{x-z}+(-1)^{k+1}\int_{\xi_{k+1}}^{\xi_{k+2}}\frac{\ud z}{x-z}+\cdots+(-1)^{n}\int_{\xi_n}^{1}\frac{\ud z}{x-z}=\\
=\frac{M_1}{x^{n+1}}+\frac{M_2}{x^{n+2}}+\cdots
\end{IEEEeqnarray*}
holds.

Thus,
\begin{equation*}
 \frac{(x+1)(x-\xi_2)^2(x-\xi_4)^2\cdots}{(x-\xi_1)^2(x-\xi_3)^2(x-\xi_5)^2\cdots}
=e^{\frac{M_1}{x^{n+1}}+\frac{M_2}{x^{n+2}}+\cdots}=1+\frac{M_1}{x^{n+1}}+\cdots
\end{equation*}
and hence, if we put
\begin{IEEEeqnarray}{c}
 (x-\xi_2)(x-\xi_4)\cdots(x-\xi_{2m})=U_{m}(x)\IEEEnonumber\\
\label{eqn:five}
(x-\xi_1)(x-\xi_3)\cdots(x-\xi_{2m+1})=V_{m+1}(x)
\end{IEEEeqnarray}
for $n=2m+1$, we obtain 
\begin{equation}
  \tag{$5_1$}
\label{eqn:fiveone}
(x+1)(x-\beta)(x-1)U_m^2(x)-(x-\alpha)V^2_{m+1}(x)=Ax+B
\end{equation}
for odd $k$ and 
\begin{equation}
 \tag{$5_2$}
\label{eqn:fivetwo}
(x+1)(x-\alpha)(x-1)U_m^2(x)-(x-\beta)V^2_{m+1}(x)=Ax+B
\end{equation}
for even $k$. On the other hand, if $n=2m$ is an even number and if we put
\begin{IEEEeqnarray}{c}
 (x-\xi_2)(x-\xi_4)\cdots(x-\xi_{2m})=P_{m}(x)\IEEEnonumber\\
\label{eqn:six}
(x-\xi_1)(x-\xi_3)\cdots(x-\xi_{2m-1})=Q_{m}(x),
\end{IEEEeqnarray}
we analogously obtain
\begin{equation}
 \tag{$6_1$}
\label{eqn:sixone}
(x+1)(x-\beta)P_m^2(x)-(x-\alpha)(x-1)Q^2_m(x)=Ax+B
\end{equation}
for even $k$ and 
\begin{equation}
 \tag{$6_2$}
\label{eqn:sixtwo}
(x+1)(x-\alpha)P_m^2(x)-(x-\beta)(x-1)Q^2_m(x)=Ax+B
\end{equation}
for odd $k$.

\S 4. We will now solve the equations from the previous paragraph by means of elliptic functions.

In doing so, we will confine ourselves to the case where $n=2m+1$.

We put
\begin{equation}
 \label{eqn:seven}
k^2=\frac{2(\beta-\alpha)}{(1-\alpha)(1+\beta)}
\end{equation}
and take $k$ ($0<k<1$) as the modulus of the Jacobi functions. Furthermore, we determine $\rho$ from the equation
\begin{equation}
 \label{eqn:eight}
\alpha=1-2\sn^2\rho
\end{equation}
under the additional condition
\begin{equation*}
 0<\rho<K.
\end{equation*}
It then follows from (\ref{eqn:seven}) and (\ref{eqn:eight}) that
\begin{equation}
 \label{eqn:nine}
\beta=2\frac{\cn^2\rho}{\dn^2\rho}-1.
\end{equation}
Finally, we set
\begin{equation}
 \label{eqn:ten}
x=\frac{\sn^2u\cdot\cn^2\rho+\cn^2u\cdot\sn^2\rho}{\sn^2u-\sn^2\rho},
\end{equation}
such that
\begin{IEEEeqnarray}{rClrCl}
 x+1 &=& \frac{2\sn^2u\cdot \cn^2\rho}{\sn^2u-\sn^2\rho},\qquad & x-\beta&=&\frac{(1-\beta^2)\dn^2u\cdot\dn^2\rho}{2(1-k^2)(\sn^2u-\sn^2\rho)},
\IEEEnonumber\\
\label{eqn:eleven}
x-\alpha &=& \frac{1-\alpha^2}{2(\sn^2u-\sn^2\rho)},\qquad &x-1&=& \frac{2\sn^2\rho\cdot\cn^2u}{\sn^2u-\sn^2\rho}.
\end{IEEEeqnarray}
The relation (\ref{eqn:ten}) is the conformal mapping from the two-sheeted Riemann surface with branch points
$D$ ($x=-1$), $A$ ($x=\alpha$), $B$ ($x=\beta$), $C$ ($x=1$) and lines of transition $DA$, $BC$ onto the period parallelogram
with vertices 
\begin{equation*}
 u=K\pm\ii K^{\prime},\quad -K\pm\ii K^{\prime}.
\end{equation*}
Thereby, the rectangle $AB^{\prime\prime}B^{\prime\prime\prime}A^{\prime}$ (Fig.~1) corresponds to the upper sheet of the Riemann
surface and the rectangle $ABB^{\prime}A^{\prime}$ to the lower one.
 \begin{figure}[!hbtp]
 \centering
\includegraphics[width=\textwidth]{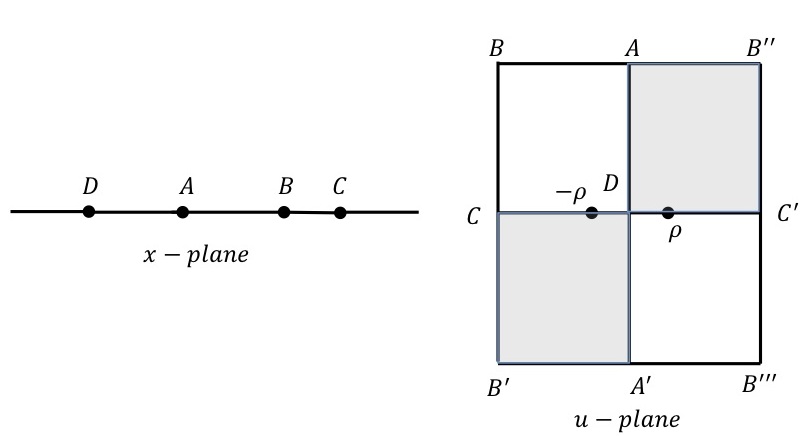}
Fig. 1
\end{figure}
We consider the function
\begin{equation*}
 \Phi(u)=\sqrt{\frac{(x+1)(x-\beta)(x-1)}{x-\alpha}}U_m(x)-V_{m+1}(x),
\end{equation*}
where $U_m(x)$, $V_{m+1}(x)$ are polynomials subject to the equation~(\ref{eqn:fiveone}), and take the square root in such a way,
that $\Phi(u)$ remains finite on the upper sheet of the Riemann surface and thus, also on the rectangle 
$AB^{\prime\prime}B^{\prime\prime\prime}A^{\prime}$.

It followos from this definition of the function $\Phi(u)$ that the point $u=\rho$ is a root of multiplicity $(m+1)$ of the function 
$\Phi(u)$.

On the other hand, the equations in (\ref{eqn:eleven}) show that
\begin{equation*}
 \Phi(u)=\frac{2\sn\rho\cdot\cn\rho\cdot\sn u\cdot\cn u\cdot \dn u}{\sn^2u-sn^2\rho}U_{m}(x)-V_{m+1}(x).
\end{equation*}
This implies that, first of all,
\begin{equation*}
 \Phi(u+2K)=\Phi(u+2\ii K^{\prime})=\Phi(u),
\end{equation*}
and, secondly,
\begin{equation*}
 \Phi(-u)=-\sqrt{\frac{(x+1)(x-\beta)(x-1)}{x-\alpha}}
U_{m}(x)-V_{m+1}(x)=-\frac{Ax+B}{(x-\alpha)\Phi(u)}.
\end{equation*}
This last property shows us that the point $u=-\rho$ is a pole of multiplicity $(m+1)$ of $\Phi(u)$; moreover, $\Phi(u)$ has a simple pole
in $u=\ii K^{\prime}$.

All these facts and the double periodicity of the function $\Phi(u)$ lead to the following representation of this function:
\begin{equation*}
 \Phi(u)= \textrm{Const}\cdot\left[\frac{H(\rho-u)}{H(\rho+u)}\right]^{m+1} \frac{\Theta(u+\overline{2m+2}\rho)}{\Theta(u)}.
\end{equation*}
Therefore, we have the two equations
\begin{align*}
  V_{m+1}(x)-\sqrt{\frac{(x+1)(x-\beta)(x-1)}{x-\alpha}}U_{m}(x)&=C \left[\frac{H(\rho-u)}{H(\rho+u)}\right]^{m+1}\frac{\Theta(u+\overline{2m+2}\rho)}{\Theta(u)}\\
  V_{m+1}(x)+\sqrt{\frac{(x+1)(x-\beta)(x-1)}{x-\alpha}}U_{m}(x)&=C \left[\frac{H(\rho+u)}{H(\rho-u)}\right]^{m+1}\frac{\Theta(u-\overline{2m+2}\rho)}{\Theta(u)}
\end{align*}
and, thus, we may present the solution to the undetermined equation (\ref{eqn:fiveone}) in the following form
\begin{equation}
\tag{$12_1$}
 \label{eqn:twelveone}
V_{m+1}(x)=\frac{C}{2}
\left\{
\left[\frac{H(\rho-u)}{H(\rho+u)}\right]^{m+1}\frac{\Theta(u+\overline{2m+2}\rho)}{\Theta(u)}+\left[\frac{H(\rho+u)}{H(\rho-u)}\right]^{m+1}
\frac{\Theta(u-\overline{2m+2}\rho)}{\Theta(u)}
\right\}
\end{equation}
\begin{IEEEeqnarray*}{rcl}
U_m(x)&=&\frac{C}{2}\sqrt{\frac{x-\alpha}{(x+1)(x-\beta)(x-1)}}
\Bigg\{
\left[\frac{H(\rho+u)}{H(\rho-u)}\right]^{m+1}\frac{\Theta(u-\overline{2m+2}\rho)}{\Theta(u)}\\
&&-\left[\frac{H(\rho-u)}{H(\rho+u)}\right]^{m+1}\frac{\Theta(u+\overline{2m+2}\rho)}{\Theta(u)}
\Bigg\}.
\end{IEEEeqnarray*}
In the same way, we obtain that the solution to the equation (\ref{eqn:fivetwo}) is of the form
\begin{equation}
\tag{$12_2$}
 \label{eqn:twelvetwo}
V_{m+1}(x)=\frac{D}{2}
\left\{
\left[\frac{H(\rho-u)}{H(\rho+u)}\right]^{m+1}\frac{\Theta_1(u+\overline{2m+2}\rho)}{\Theta_1(u)}+\left[\frac{H(\rho+u)}{H(\rho-u)}\right]^{m+1}
\frac{\Theta_1(u-\overline{2m+2}\rho)}{\Theta_1(u)}
\right\}
\end{equation}
\begin{IEEEeqnarray*}{rcl}
U_m(x)&=&\frac{D}{2}\sqrt{\frac{x-\beta}{(x+1)(x-\alpha)(x-1)}}
\Bigg\{
\left[\frac{H(\rho+u)}{H(\rho-u)}\right]^{m+1}\frac{\Theta_1(u-\overline{2m+2}\rho)}{\Theta_1(u)}\\
&&-\left[\frac{H(\rho-u)}{H(\rho+u)}\right]^{m+1}\frac{\Theta_1(u+\overline{2m+2}\rho)}{\Theta_1(u)}
\Bigg\}.
\end{IEEEeqnarray*}

It is  evident that, on the other hand, due to the transformation (\ref{eqn:ten}), the functions in $u$ on the right hand-side in the above formulas
turn out to be polynomials in $x$, in fact, polynomials of degrees $m+1$ and $m$, respectively. However, a priori, it is not known
that the zeros of these polynomials satisfy the inequalities discussed in \S~3.

On the basis of our foregoing considerations, however, we can already at this point claim that, in any case, at least one
pair of functions (\ref{eqn:twelveone}) or (\ref{eqn:twelvetwo}) contribute to the solution 
\begin{equation*}
 f(x)=U_m(x)V_{m+1}(x)
\end{equation*}
to our problem for $n=2m+1$.

In the following section we will show how one can determine which pair of functions (\ref{eqn:twelveone}), (\ref{eqn:twelvetwo})
solves our problem, given $\alpha$ and $\beta$.

\S 5. Let 
\begin{equation*}
 (2m+2)\rho=pK+\sigma,
\end{equation*}
where $p$ is an integer and $\sigma$ is subject to the inequality
\begin{equation*}
 0\leq\sigma\leq K.
\end{equation*}

It is expedient to distinguish between the following three cases: 1) $0<\sigma<K$, $p$ --- odd number; 2)  $0<\sigma<K$, $p$ --- even number;
3) $\sigma=0$ or $\sigma=K$.

\textit{First case:} $0<\sigma<K$, $p$ --- odd number.

We take the functions (\ref{eqn:twelvetwo}) and rewrite them in the form
\begin{IEEEeqnarray*}{c}
 V_{m+1}(x)=\frac{D}{2}
\left[\frac{H(\rho+u)}{H(\rho-u)}\right]^{m+1}\frac{\Theta_1(u-\overline{2m+2}\rho)}{\Theta_1(u)}
\left\{1+\Omega(u) \right\},\\
U_m(x)=\frac{D}{2}\sqrt{\frac{x-\beta}{(x+1)(x-\alpha)(x-1)}}
\left[\frac{H(\rho+u)}{H(\rho-u)}\right]^{m+1}\frac{\Theta_1(u-\overline{2m+2}\rho)}{\Theta_1(u)}\left\{1-\Omega(u) \right\},
\end{IEEEeqnarray*}
where we put
\begin{equation*}
 \Omega(u)= \left[\frac{H(\rho-u)}{H(\rho+u)}\right]^{2m+2}\frac{\Theta_1(u+\overline{2m+2}\rho)}{\Theta_1(u-\overline{2m+2}\rho)}.
\end{equation*}
Simple calculations show that
\begin{equation*}
 \Omega(0)=1,\quad \Omega(\ii K^{\prime})=1,\quad \Omega(K+\ii K^{\prime})=-1,\quad \Omega(K)=1;
\end{equation*}
 however, one can  easily check that $x=-1,~\alpha,~\beta,~1$ are not zeros of the polynomials $V_{m+1}(x)$, $U_{m}(x)$.

Moreover, we see that, in the intervals $[-1,\alpha]$, $[\beta,1]$, the function $\Omega(u)$ equals one in modulus.

The zeros of the polynomials $V_{m+1}(x)$, $U_{m}(x)$ in the intervals $(-1,\alpha)$, $(\beta,1)$ thus correspond to those points
on the open line segments $DA$ and $B^{\prime\prime}C^{\prime}$, where $\Omega(u)=\pm1$.

In order to find the number of these points, we will consider the path $DAMB^{\prime\prime}C^{\prime}ND$ (Fig.~2); and we set $\omega=\ii K^{\prime}+\sigma$. By running through this path $\arg\Omega(u)$ increases by
\begin{equation*}
 -(2m+2)2\pi+2\pi.
\end{equation*}

Putting
\begin{equation*}
 \arg\Omega(0)=0
\end{equation*}
hence gives
\begin{equation*}
 \arg\Omega(K)=-\pi(2m+2)+2\pi=-2m\pi.
\end{equation*}
\begin{figure*}[hbtp]
 \centering
\includegraphics[width=0.7\textwidth]{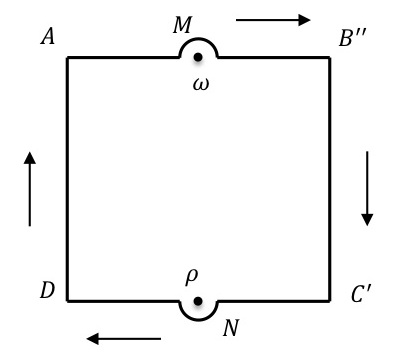}\\
Fig. 2
\end{figure*}
Furthermore, we have
\begin{equation*}
 \arg\Omega(\ii K^{\prime}+K)=\arg\Omega(\ii K^{\prime})+\pi
\end{equation*}
and
\begin{equation*}
 \arg\Omega(\ii K^{\prime})=(2m+2)\arg\frac{H(\rho-u)}{H(\rho+u)}+\arg\left.\frac{\Theta(u+\sigma)}{\Theta(u-\sigma)}\right|_{u=\ii K^{\prime}}.
\end{equation*}
But on the other hand we have
\begin{align*}
 \frac{H(\rho-\ii K^{\prime})}{H(\rho+\ii K^{\prime})}&=-e^{\frac{\pi\ii}{K}\rho}\\
\frac{\Theta(\ii K^{\prime}+\sigma)}{\Theta(\ii K^{\prime}-\sigma)}&=-e^{-\frac{\pi\ii}{K}\sigma},
\end{align*}
and, since $\arg\frac{H(\rho-u)}{H(\rho+u)}$ has to be negative and smaller than $2\pi$ in modulus in the point $A$ when approaching along
$DA$, we have
\begin{equation*}
 \arg\left.\frac{H(\rho-u)}{H(\rho+u)}\right|_{u=\ii K^{\prime}}=-\pi+\frac{\pi}{K}\rho;
\end{equation*}
similarly, one finds that
\begin{equation*}
 \arg\left.\frac{\Theta(u+\sigma)}{\Theta(u-\sigma)}\right|_{u=\ii K^{\prime}}= \pi-\frac{\pi}{K}\sigma.
\end{equation*}
Thus,
\begin{equation*}
 \arg\Omega(\ii K^{\prime})=-(2m+2)\pi+\pi+\frac{\pi}{K}\left[(2m+2)\rho-\sigma\right]=-\pi(2m+1-p),
\end{equation*}
and hence
\begin{equation*}
 \arg\Omega(\ii K^{\prime}+K)=-\pi(2m-p).
\end{equation*}
The change that $\arg\Omega(u)$ undergoes along $DA$ is $\pi(2m+1-p)$, and along the line $B^{\prime\prime}C^{\prime}$ this
change amounts to $\pi p$.

Now, $p$ is an odd number, let $p=2q+1$; then each of the polynomials $U_m(x)$, $V_{m+1(x)}$ have $q$ zeros in the interval $(\beta,1)$
and these zeros are simple: if there were further zeros of these polynomials  in this
interval, then we would either have $\Omega(u_0)=\pm1$ and $\Omega^{\prime}(u_0)=0$ for some point ($u=u_0$) in this interval, or the
equations $\Omega(u_1)=\Omega(u_2)=\pm1$ would have to hold for two points ($u=u_1$, $u=u_2$), while $\Omega(u)$ would differ from $+1$
between them; in the first case, the polynomial would have a multiple zero, and, in the second case, we would have a simple zero of one
of the polynomials without a zero of the other in between; both possibilities are in contradiction to the properties of the zeros of the
polynomial $f(x)$ we are looking for.

Similarly, we see that the number of roots from the second interval, which may satisfy the inequalities from \S 3, are equal to
$m-q-1$, $m-q$. Therefore, in both intervals, we get $m-1$ zeros for one of the polynomials and $m$ zeros for the other, which satisfy
the inequalities from \S~3. The sum is smaller than $2m+1$.

These considerations show that, in this case, the polynomials (\ref{eqn:twelvetwo}) cannot define a solution, consequently, the solution
has to equal the product of the polynomials (\ref{eqn:twelveone}), and, indeed, similar considerations show that the polynomials
$U_m(x)$ and $V_{m+1}(x)$ each have $m-q$ zeros in the interval $(-1,\alpha)$, whereas the respective numbers are $q$, $q+1$ in $(\beta,1)$,
while the inequalities from \S~3 are satisfied.

\textit{Second case: } $0<\sigma<K$ --- $p$ even number.

By similar considerations we find that, in this case, the solution to the problem equals the product of the polynomials (\ref{eqn:twelvetwo}),
and, if $p=2q$, then the solution to the problem has
\begin{equation*}
 2m-2q+1
\end{equation*}
zeros in the interval $(-1,\alpha)$.

The polynomials (\ref{eqn:twelveone}) do not give a solution in this case.

\textit{Third case: } $(2m+2)\rho\equiv0\pmod K$.

Here, we have a limiting case of the first as well as of the second case.

Thus, we obtain two solutions:
\begin{align*}
 f_1(x)&=\sqrt{\frac{x-\alpha}{(x+1)(x-\beta)(x-1)}}
\left\{ 
\left[ \frac{H(\rho-u)}{H(\rho+u)} \right]^{2m+2}- \left[ \frac{H(\rho+u)}{H(\rho-u)} \right]^{2m+2}
\right\},\\
f_2(x)&=\sqrt{\frac{x-\beta}{(x+1)(x-\alpha)(x-1)}}
\left\{ 
\left[ \frac{H(\rho-u)}{H(\rho+u)} \right]^{2m+2}- \left[ \frac{H(\rho+u)}{H(\rho-u)} \right]^{2m+2}
\right\}.
\end{align*}

These solutions are of the form
\begin{align*}
 f_1(x)&=(x-\alpha)\varphi(x),\\
f_2(x)&=(x-\beta)\varphi(x),
\end{align*}
where $\varphi(x)$ is a polynomial of degree $n-1$, and they are different; as one can easily see, our problem has infinitely many
solutions
\begin{equation*}
 (x-\vartheta)\varphi(x)\qquad(\alpha\leq\vartheta\leq\beta)
\end{equation*}
in the case in question.

From these solutions we pick the following
\begin{equation*}
 f_0(x)= A
\frac{x-\gamma}{\sqrt{(x+1)(x-\alpha)(x-\beta)(x-1)}}
\left\{ 
\left[ \frac{H(\rho-u)}{H(\rho+u)} \right]^{2m+2}- \left[ \frac{H(\rho+u)}{H(\rho-u)} \right]^{2m+2}
\right\},
\end{equation*}
where
\begin{equation*}
 \gamma=\alpha+\frac{2\sn\rho\cdot\cn\rho}{\dn\rho}\frac{\Theta^{\prime}(\rho)}{\Theta(\rho)}
\end{equation*}
(this number is between $\alpha$ and $\beta$, see \cite{AchUeb1933}).

The significance of this special choice lies in the fact that (cf~\cite{AchUeb1933})
\begin{equation*}
 \frac{2(x-\gamma)}{\sqrt{(1-\alpha)(1+\beta)}} = \frac{H^{\prime}(\rho-u)}{H(\rho-u)}-\frac{H^{\prime}(\rho+u)}{H(\rho+u)}.
\end{equation*}
Consequently,
\begin{equation}
\addtocounter{equation}{1}
\label{eqn:thirteen}
 f_0(x)=B\frac{\ud}{\ud x}
\left\{ 
\left[ \frac{H(\rho-u)}{H(\rho+u)} \right]^{2m+2}- \left[ \frac{H(\rho+u)}{H(\rho-u)} \right]^{2m+2}
\right\},
\end{equation}
and, since
\begin{equation*}
 (2m+2)\rho=pK,
\end{equation*}
the above polynomial inside the braces differs only by a constant factor from the polynomial
\begin{equation*}
 T_{2m+2}(x;p,k),
\end{equation*}
which I investigated in my earlier papers (\cite{AchUeb1933}). This polynomial was similar to the 
Tch\' ebyschef polynomial $T_n(x)=\cos n\arccos x$ in many ways; now, we have found that $T_n(x)$ shares one further common property.

R~e~m~a~r~k.

Easy calculations show that the coefficient $B$ from formula~(13)\footnote{ 
The entity
\begin{equation*}
 \tau=\frac{1}{2} \left[ \frac{\Theta(0)\Theta_1(0)}{\Theta(\rho)\Theta_1(\rho)} \right]^{2}
\end{equation*}
is the transfinite diameter (M. Fek\' ete) of the point set $E$ (cf \cite{AchUeb1933}).

In the special case where
\begin{equation*}
 \alpha=-\beta
\end{equation*}
we have 
\begin{equation*}
 \tau=\frac{\sqrt{1-\alpha^2}}{2}.
\end{equation*}
}
 has the value
\begin{equation}
 \tag{$13_1$}
\label{eqn:thirteenone}
B=\frac{1}{(2m+2)4^{m+1}}\left[ \frac{\Theta(0)\Theta_1(0)}{\Theta(\rho)\Theta_1(\rho)} \right]^{4m+4}=
\frac{\tau^{2m+2}}{2m+2},
\end{equation}
moreover,
\begin{align}
 \int_{E}|f_0(x)|\ud x&=
B\int_{x\in E}\left| \ud \left\{ \left[ \frac{H(\rho-u)}{H(\rho+u)} \right]^{2m+2} + \left[ \frac{H(\rho+u)}{H(\rho-u)} \right]^{2m+2}  \right\}  \right|
\nonumber\\
&=2B \int_{0}^{\pi}| \ud\cos(2m+2) v|= \int_{0}^{(2m+2)\pi}|\ud\cos v|\nonumber\\
\tag{$13_2$} \label{eqn:thirteentwo}
&=2B\cdot(2m+2)\cdot2\cdot\int_{0}^{\frac{\pi}{2}}|\ud\cos u|=4B(2m+2).
\end{align}
It now follows from the Fundamental Theorem on polynomials with minimal deviation that $f_0(x)$ deviates least from zero in $E$
among all polynomials of the form $x^{2m+1}+p_1x^{2m}+\cdots$
w.r.t. the weight
\begin{equation*}
 \frac{ \sqrt{(1-x^2)(\alpha-x)(\beta-x)} }{|x-\gamma|}
\end{equation*}
in the sense of Tsch\' ebyschef,
and the deviation attains the value
\begin{equation*}
 2B(2m+2).
\end{equation*}
Thus, for $(2m+2)\rho\equiv0\pmod K$, the exact relation of S.~N.~Bernstein
\begin{IEEEeqnarray*}{c}
 \min\int_{E}\left| x^{2m+1} + p_1x^{2m} +\cdots+p_{2m+1} \right|\ud x=\\
2\min\left\{  \max_{x\in E}\left| \left( x^{2m+1} + p_1x^{2m} +\cdots+p_{2m+1} \right) 
\frac{ \sqrt{(1-x^2)(\alpha-x)(\beta-x)} }{x-\gamma}
  \right| \right\}=\\
=2\min\left\{ \max_{x\in E}\left|x^{2m+2} + q_1x^{2m+1} +\cdots+q_{2m+2}   \right| \right\}
\end{IEEEeqnarray*}
holds.

\S 6. In the previous paragraph we have fully solved our problem for $n=2m+1$.

This result can also be expressed differently. To this end, we introduce some notation. Let $\alpha$ be fixed whereas we consider $K$
as a function of the variable argument $k$; assume that the equation
\begin{equation}
\label{eqn:fourteen}
 \alpha=1-2\sn^2\left(\frac{p}{2m+2}K,k\right)
\end{equation}
is solved w.r.t. $k$, where $p$ denotes any of the numbers $1,2,\ldots,q$ and $q<2m+2$ is the largest natural number such that
(\ref{eqn:fourteen}) still leads to a positive fraction $k$ (so, if $\cos\frac{\pi\mu}{2m+2}\leq\alpha$, then $q$ cannot be larger
than $\mu$). 

We therefore obtain the quantities $k_1,~k_2,\ldots,~k_q$ and, subsequently, we find the corresponding values for $\beta$
\begin{equation*}
 \beta_q<\beta_{q-1}<\ldots<\beta_1<\beta_0=1
\end{equation*}
through formula (\ref{eqn:seven}).

Setting $\beta_{q+1}=\alpha$ we may then claim that our problem has a unique solution for
\begin{equation*}
 \beta_{2j}<\beta<\beta_{2j-1}\qquad(j=1,2,\ldots;~2j\leq q+1),
\end{equation*}
which is defined by the formulas in (\ref{eqn:twelveone}); for 
\begin{equation*}
  \beta_{2j+1}<\beta<\beta_{2j}\qquad(j=0,1,\ldots;~2j+1\leq q+1)
\end{equation*}
too there exists only one solution, but it is defined by the formulas in (\ref{eqn:twelvetwo}); and only for 
\begin{equation*}
 \beta=\beta_j
\end{equation*}
there are infinitely many solutions, among which there is one related to ``elliptic'' polynomials as well, like the solution to the 
Korkine-Zolotaref problem with the trigonometric polynomial $\cos n\arccos x$ does.

It is not hard to give asymptotical estimates (for $m\to\infty$) of
\begin{equation*}
 \min\int_{E}\left|  x^{2m+1}+p_1x^{2m}+\cdots+p_{2m+1} \right|\ud x,
\end{equation*}
in case $\alpha$ and $\beta$ are fixed numbers.

For any $m$ one can find two numbers $\beta^{(m)}_{j+1}$, $\beta^{(m)}_{j}$ such that
\begin{equation*}
 \beta^{(m)}_{j+1}\leq\beta\leq\beta^{(m)}_{j},
\end{equation*}
and where
\begin{equation*}
 \lim_{m\to\infty}\beta^{(m)}_{j}=\lim_{m\to\infty}\beta^{(m)}_{j+1}=\beta\qquad(j=j_{(m)}).
\end{equation*}
Now, if 
\begin{equation*}
 G_{2m+1}(b)=\min\left\{ \int_{-1}^{\alpha}\left|  x^{2m+1}+p_1x^{2m}+\cdots \right|\ud x
+\int_{b}^{1}\left|  x^{2m+1}+p_1x^{2m}+\cdots \right|\ud x \right\},
\end{equation*}
we obviously have that
\begin{IEEEeqnarray*}{rCl}
 G_{2m+1}(\beta) &\geq&  G_{2m+1}(\beta^{(m)}_{j})\\
 G_{2m+1}(\beta^{(m)}_{j+1}) &\geq&  G_{2m+1}(\beta),
\end{IEEEeqnarray*}
and, thus, the following inequality holds
\begin{IEEEeqnarray*}{rcl}
 \IEEEeqnarraymulticol{3}{l}{
\left[ \frac{\Theta(0;k_j)\Theta_1(0;k_j)}{\Theta\left(\frac{j}{2m+2}K_j;k_j\right)\Theta_1\left(\frac{j}{2m+2}K_j;k_j\right)}
 \right]^{4(m+1)}
\leq 2^{2m}G_{2m+1}(\beta)\leq
}\\
\quad
&\leq& 
\left[ \frac{\Theta(0;k_{j+1})\Theta_1(0;k_{j+1})}{\Theta\left(\frac{j+1}{2m+2}K_{j+1};k_{j+1}\right)\Theta_1\left(\frac{j+1}{2m+2}
K_{j+1};k_{j+1}\right)}
 \right]^{4(m+1)}.
\end{IEEEeqnarray*}
From this it evidently follows that, for a sequence
\begin{equation*}
 m_1,~m_2,\ldots,
\end{equation*}
for which $\sigma\to\infty$, we have the following asymptotic equality
\begin{equation*}
 G_{2m+1}(\beta)\sim\frac{1}{2^{2m}}\left[ \frac{\Theta(0)\Theta_1(0)}{\Theta(\rho)\Theta_1(\rho)} \right]^{4(m+1)}.
\end{equation*}

\S 7. What still remains is to provide formulas for even $n$ $(=2m)$.

Here, the solution to the problem, again, is either of the form
\begin{IEEEeqnarray*}{rCl}
 f(x)&=&A\sqrt{\frac{x-\beta}{(x+1)(x-\alpha)(x-1)}}
\Bigg\{ 
\left[ \frac{H(\rho-u)}{H(\rho+u)} \right]^{2m+1}
\frac{\Theta_1^2(u+\overline{2m+1}\rho)}{\Theta_1^2(u)}
-\\
\qquad 
 &&
  - \left[ \frac{H(\rho+u)}{H(\rho-u)} \right]^{2m+1} \frac{\Theta_1^2(u-\overline{2m+1}\rho)}{\Theta_1^2(u)}
\Bigg\}
\end{IEEEeqnarray*}
or
\begin{IEEEeqnarray*}{rCl}
 f(x)&=&B\sqrt{\frac{x-\alpha}{(x+1)(x-\beta)(x-1)}}
\Bigg\{ 
\left[ \frac{H(\rho-u)}{H(\rho+u)} \right]^{2m+1}
\frac{\Theta^2(u+\overline{2m+1}\rho)}{\Theta^2(u)}
-\\
\qquad 
 &&
  - \left[ \frac{H(\rho+u)}{H(\rho-u)} \right]^{2m+1} \frac{\Theta^2(u-\overline{2m+1}\rho)}{\Theta^2(u)}
\Bigg\},
\end{IEEEeqnarray*}
in fact, it is the first one, if the integer $p$ satisfying the relation
\begin{equation*}
 (2m+1)\rho=pK+\sigma\qquad(0\leq\sigma\leq K)
\end{equation*}
is even, and the second one, if it is odd.

\newpage

\noindent
{\bf \Large Small Commentary}

\bigskip

I think this paper is extremely interesting because the famous Akhiezer's orthogonal polynomials, see  \cite[Chapter X]{Aefr} or its 
English translation \cite{Aef},   were presented here for  the first time.  At least this paper is definitely not well-known,
 in particular, because it was published in {\em Communications of the Kharkov\footnote{At that moment, Kharkov  was a city in the 
former Soviet Union, currently Ukraine. Normally, people speak Russian and/or Ukrainian.} Mathematical Society} \cite{AKH36} in 
{\em German}!

Although it is not mentioned in this paper that the polynomials $U_m(x)$ and $V_{m+1}(x)$ are orthogonal with respect to the special 
weights, we do not need to make any explanations since, in Chapter X \cite{Aefr, Aef}, they are presented and described as orthogonal
 polynomials under  the names $Q_n(x)$ and $P_n(x)$ respectively, satisfying the same functional equation \eqref{eqn:fiveone} and 
given by the same expressions in terms of theta-functions \eqref{eqn:twelveone}.

Let me mention that it was quite a non-trivial task to find this paper if one follows the footnote on page 222 in the Russian version 
of {\em Elliptic functions} \cite{Aefr}: there is not a single reference to the term {\em orthogonal polynomials}, neither in the title nor in the 
abstract of the paper written in Ukrainian (and as we can see now, even in the text of this paper it is not mentioned that the 
polynomials $U_m(x)$ and $V_{m+1}(x)$ are orthogonal with respect to the special weights! See the above paragraph.). But the worst thing is,
 that  the year of publication (actually the only hint in this footnote related to the paper) is given incorrectly, 1938 instead of 1936,
 see the front page of the volume in Fig. 3. An English-reading investigator would have had a much better fate, since \cite{Aef}, 
in contrast to \cite{Aefr}, contains a normal list of references with the correct bibliographic information.

The current translation from German  to English was given by Florian Puchhammer. Many thanks to him! 
%Many thanks to G. Feldman, who convinced me that the paper should exist, as soon as N. I. Akhiezer 
%mentioned that it exists in the footnote, even with a possible misprint in a highly incomplete bibliographic 
%information.

\bigskip

\hskip 8cm Peter Yuditskii

\begin{figure}[!hbtp]
 \centering
\includegraphics[scale=0.33]{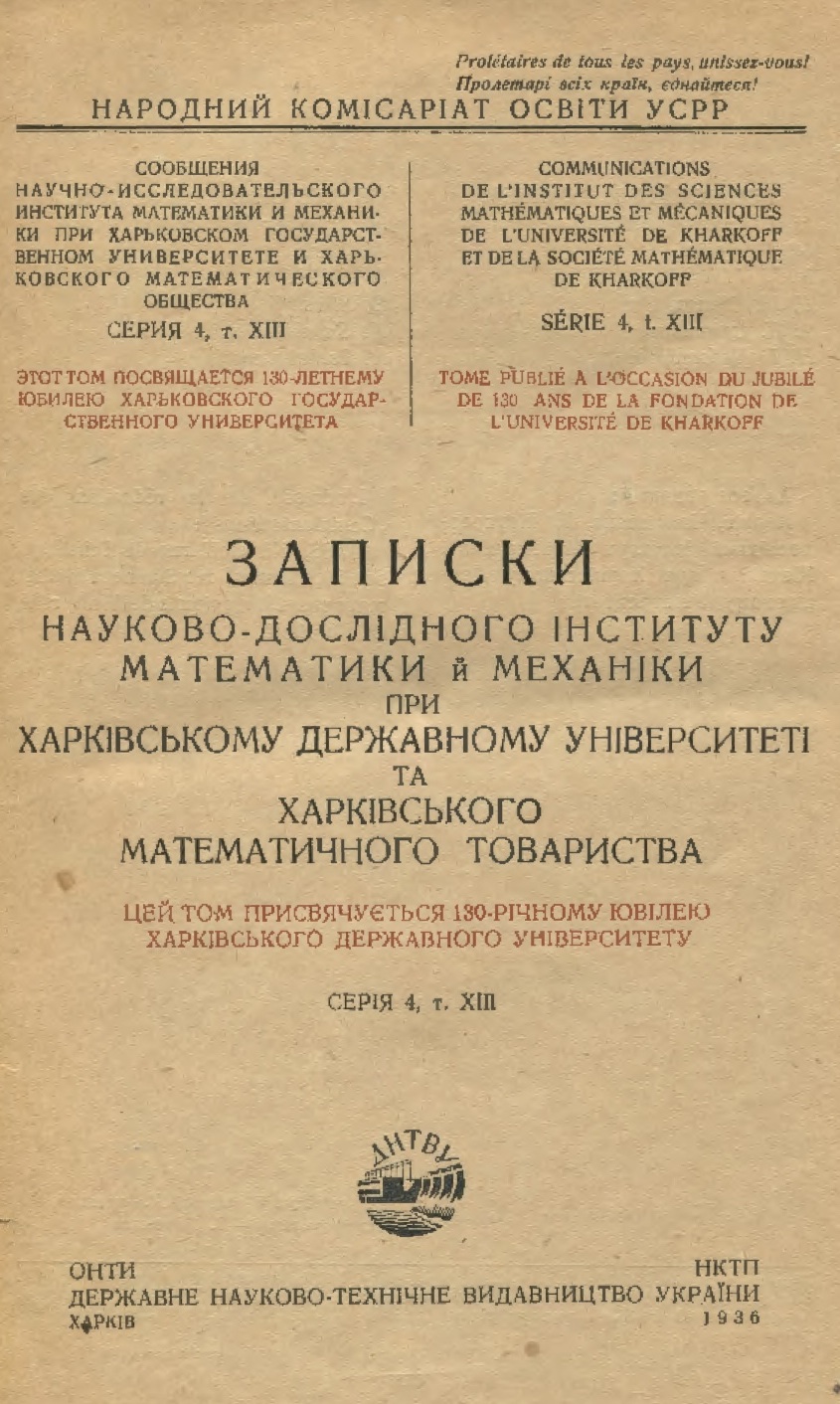}

Fig. 3
\end{figure}

\end{document}